\documentclass[12pt]{amsart}
\usepackage{amssymb}
\usepackage{amscd}
\usepackage{amsfonts}
\usepackage{latexsym}
\usepackage[all]{xy}
\usepackage{verbatim}

\newtheorem{theorem}{Theorem}[section]

\newtheorem{proposition}[theorem]{Proposition}

\theoremstyle{definition}

\newtheorem{example}[theorem]{Example}

\newtheorem{problem}[theorem]{Problem}
\newtheorem{remark}[theorem]{Remark}

\begin{document}

\title{On Hele-Shaw problems arising as scaling limits}

\author{Pavel Etingof}
\address{Department of Mathematics, Massachusetts Institute of Technology,
Cambridge, MA 02139, USA}
\email{etingof@math.mit.edu}
\maketitle

\centerline{\bf To the memory of Vladimir Markovich Entov}

\section{Introduction}

This work is motivated by the paper \cite{LP}. 
In this paper, the authors consider the scaling limits 
of three discrete aggregation models with multiple sources - the internal DLA
(diffusion-limited aggregation), the rotor-router, and the divisible sandpile model. 
They show that for all three models, the scaling limit is the same, and 
it is a solution of the Hele-Shaw injection problem with multiple sources.  
In two dimensions, the latter problem can be solved explicitly using the theory
of quadrature domains, which gives an explicit formula for the scaling limit. 

They also consider the same aggregation models in two dimensions 
with just one source, but with an additional condition on the positive (horizontal) half-axis
(namely, the condition that a particle hitting the positive half-axis is killed, 
or the condition that it is directed downward). In these cases, the 
existence of the scaling limit remains unknown, although it is expected to exist, 
and computer-generated pictures for its shape are given in Fig. 4 of \cite{LP}. 

The goal of this paper is to study the Hele-Shaw problems which are expected to 
be the scaling limits of the above models with a condition. Namely, we 
provide an explicit solution for the Hele-Shaw problem corresponding to the killing
condition, which is a close fit with the left shape at Fig. 4 of \cite{LP}. 
We also describe moment properties of the solution of the Hele-Shaw problem corresponding 
to the downward condition (the right shape at Fig. 4 of \cite{LP}), although 
we are unable to compute this shape explicitly. 

{\bf Acknowledgements.} This paper is dedicated to the memory 
of my teacher Vladimir Markovich Entov, who introduced me to the 
subject of Hele-Shaw flows. He was an extraordinary person and scientist, 
and interaction with him was one of my best experiences.

I am very grateful to Lionel Levine for introducing me to the problem, 
and for checking that my formulas fit the results of computer simulations. 
Without discussions with him, this paper would not have been written. 

This work was partially supported by the NSF grant DMS-9988796.

\section{Flows with the killing/reflecting condition}

\subsection{The killing condition on a half-axis}

Consider a 2-dimensional discrete aggregation model with injection at the origin and killing condition on a half-axis; 
we prefer to impose this condition on the negative half-axis (rather than the positive half-axis used in \cite{LP}). 
To be more concrete, consider the internal DLA model. In this model, particles are emitted from the origin 
into a rectangular lattice, and each particle walks randomly on this lattice until it reaches an unoccupied node, 
where is stays, and a new particle is emitted. In addition, one stipulates that a particle that hits the negative half-axis 
is killed (after which a new particle is emitted). The scaling limit of this model is the limit in which the number $N$ of surviving particles
\footnote{Note that most of the particles will not survive. Namely, as was pointed out to us by Lionel Levine, 
it follows from the Beurling estimate that for $N$ particles to survive, the number of produced particles 
will have to be $\sim C N^{5/4}$, where $C$ is a constant.}
goes to infinity, and the spacing of the lattice goes to zero as $(A/N)^{1/2}$, where $A$ is a constant. If the scaling limit exists, it is a region 
$\Omega$ in the plane of area $A$, containing the origin in its interior. 
We would like to compute this region explicitly using conformal mappings, assuming the scaling limit
exists and is the solution of the Hele-Shaw problem naturally associated to this model.    

Now consider the Hele-Shaw problem for the asymptotic region $\Omega$.
We refer the reader to the books \cite{VE} and \cite{GV}
for the basics on Hele-Shaw moving boundary problems.  

\begin{problem}\label{pr1} Let $\Phi$ be the harmonic potential in $\Omega$ with deleted
negative half-axis, with boundary conditions $\Phi=0$ on the boundary
$\Gamma$ of $\Omega$ and also on the negative half-axis (which is the continuous 
counterpart of the killing condition on the negative half-axis in the discrete model), 
and $\Phi\sim -{\rm Re}(z^{-1/2})$ near zero (where $z$ is the complex coordinate in the plane). 
The condition on $\Omega$ is that the infinitesimal deformation of the boundary given by the normal
derivative of $\Phi$ is homothetic. 
\end{problem}

Let us find $\Omega$ from this data (up to scaling). For this 
purpose, let $f: D\to \Omega$ be the conformal map of the unit
disk onto $\Omega$ with $f(0)=0$ and $f'(0)>0$ (it is unique). 
Let $\Psi(u)=\Phi(f(u))$. Then $\Psi$ is a harmonic potential 
on the disk with a cut along the negative half-axis, and 
$\Psi$ behaves like $-{\rm Re}(u^{-1/2})$ at zero. 
Thus, up to scaling, 
$$
\Psi={\rm Re}(u^{1/2}-u^{-1/2}).
$$
(we choose the standard branch of the square root with a cut on
the negative half-axis). Indeed, $\Psi$ is uniquely determined 
by its properties, and the given $\Psi$ satisfies them. 

Now, if $|\zeta|=1$ then 
the velocity vector of the boundary at the point $f(\zeta)$ 
under the transformation by the normal derivative of $\Phi$ is 
$\frac{\partial \Psi}{\partial n}|f'(\zeta)|^{-1}$, 
while the outward unit normal to the boundary at $f(\zeta)$ 
is $\frac{\zeta f'(\zeta)}{|f'(\zeta)|}$. Thus 
the condition of homothetic deformation has the form 
$$
{\rm Re}\left(\frac{\zeta f'(\zeta)}{|f'(\zeta)|}\overline{f(\zeta)}\right)=
\frac{\partial \Psi}{\partial n}|f'(\zeta)|^{-1}, |\zeta|=1,
$$
i.e., 
\begin{equation}
{\rm Re}(\zeta f'(\zeta)\overline{f(\zeta)})=
\frac{\partial \Psi}{\partial n}, |\zeta|=1,
\end{equation}
or 
\begin{equation}\label{maineq}
{\rm Re}(\zeta f'(\zeta)\overline{f(\zeta)}) 
=\frac{1}{2}(\zeta^{1/2}+\zeta^{-1/2}), |\zeta|=1.
\end{equation}

Now we would like to solve (\ref{maineq}).
We use the method similar to one described in 
\cite{EEK}, for evolution of polygons (Section 9). 

Let $D=\zeta\partial_\zeta$. 
Also for any function $g$ set $g^*(\zeta)=g(1/\zeta)$, so 
$(Dg)^*=-D(g^*)$. 

Note that $f$ is real, i.e. it has real Taylor coefficients at
$0$. Thus equation (\ref{maineq})
can be rewritten in the form 
\begin{equation}\label{maineq1}
Df\cdot f^*+f\cdot (Df)^*=\zeta^{1/2}+\zeta^{-1/2}, |\zeta|=1. 
\end{equation}
Differentiating this, we get 
\begin{equation}\label{maineq2}
D^2f\cdot f^*-f\cdot (D^2f)^*=\frac{1}{2}(\zeta^{1/2}-\zeta^{-1/2}).
\end{equation}

Let 
$$
h=\frac{D^2f-\frac{1}{2}\frac{\zeta^{1/2}-\zeta^{-1/2}}
{\zeta^{1/2}+\zeta^{-1/2}}Df}{f}.
$$
Then equations (\ref{maineq1},\ref{maineq2}) imply that 
$$
h-h^*=0, |\zeta|=1. 
$$
But by its definition, $h$ extends to a function which is 
holomorphic everywhere in the unit disk except possibly the
point $-1$. Hence $h$ is holomorphic in the Riemann sphere with
the possible exception of the point $-1$. Also, it's clear that 
$Df(-1)=0$, which by the removable singularity theorem 
and Liouville's theorem implies that in fact $h$ is a constant. 

Thus, we find that $f$ satisfies the differential equation 
$$
D^2f+\frac{1}{2}\frac{1-\zeta}
{1+\zeta}Df=hf
$$
The constant $h$ is easy to find from the condition that 
$f(\zeta)$ is proportional to $\zeta+O(\zeta^2)$: 
namely, $h=\frac{3}{2}$. 

From this is it easy to solve the differential equation 
by the power series method. The answer is, up to scaling 
$$
f(\zeta)=\frac{15}{16}\sum_{k\ge 1}\frac{k}{(k^2-1/4)(k^2-9/4)}(-\zeta)^k.
$$
(the normalization is such that $f'(0)=1$). Thus, we have 

\begin{proposition}\label{prop1}
The region solving Problem \ref{pr1} is defined by the conformal mapping 
\begin{equation}\label{negax}
f(\zeta)=\frac{15}{32}(1+\zeta)^2(\zeta^{-1}-
(1-\zeta)\zeta^{-3/2}{\rm arctan}\zeta^{1/2})-\frac{5}{8}.
\end{equation}
\end{proposition}

Thus we have a ``logarithmic cusp'' near $\zeta=-1$,
i.e. $f(-1+u)$ behaves like $cu^2\log u$ (as opposed to the
ordinary semicubic cusp, with local behavior $cu^2$).  
We also see that the region is twice as thick in the positive
direction as it is in the negative direction. 

Formula (\ref{negax}) is a very close fit with the left shape on Fig. 4 of \cite{LP} 
(when it is rotated by $180^o$ to match the positive and the negative half-axes). 
This leads to a conjecture that the region defined by 
formula (\ref{negax}) is the scaling limit of the model. 

\subsection{The killing condition on the sides of an angle}

The analysis of the previous subsection can be generalized to the following more general
problem. Let $0<b\le 1$, and consider the angle $\lbrace{ z\in \Bbb C|{\rm arg}z\in
[-\pi b,\pi b]\rbrace}$. Consider the internal DLA model as above inside the angle
with a killing condition on the sides of the angle. Then the setting 
of the previous subsection is the special case $b=1$. 

The Hele-Shaw problem for the asymptotic region $\Omega$ 
in this case is as follows. 

\begin{problem}\label{pr2} Suppose that $\Omega$ is a wedge-shaped region bounded by 
the sides of the angle and some curve $\Gamma$ connecting them 
(symmetric with respect to the horizontal axis). Let $\Phi$ be the
harmonic potential in $\Omega$ which vanishes on the boundary $\Gamma$ and the
sides of the angle, and behaves like $-{\rm Re}(z^{-1/2b})$ near
zero. Then the defining condition for $\Omega$ is that $\Gamma$ transforms
homothetically under the infinitesimal deformation by the normal
derivative of $\Phi$. 
\end{problem}

\begin{remark}
Note that the condition $b\le 1$ is not essential, as for $b>1$
we can put the angle (which is $>2\pi$) on the helical
Riemann surface, which is the universal covering of the complex
plane punctured at zero. 
\end{remark}

The region $\Omega$ for any $b$ can be found similarly to the case
$b=1$ considered above. Namely, let $f$ be the conformal map 
of the circular sector defined by the inequalities 
$|z|\le 1$, $-\pi b\le {\rm arg}z\le \pi b$ (lying on the Riemann surface
if $b>1$) onto $\Omega$, such that $f$ bijectively maps straight sides to straight
sides (so $f(0)=0$), and $f'(0)>0$. Then we have an expansion 
$$
f(\zeta)=\sum_{k\ge 0}a_k \zeta^{1+k/b},
$$
for some $a_k\in \Bbb R$. 

By a method analogous to the case $b=1$, we derive a differential
equation for $f$, which has the form 
$$
D^2f+\frac{1}{2b}\frac{1-\zeta^{1/b}}
{1+\zeta^{1/b}}Df=(1+1/2b)f.
$$
Then we can solve for $f$ by the power series method, and up to
scaling we get 

\begin{proposition}\label{prop2}
The region solving Problem \ref{pr2} is defined by the conformal mapping 
\begin{equation}\label{angl}
f(\zeta)=\zeta F(\alpha,\beta,\gamma;-\zeta^{1/b}),
\end{equation}
where $F$ is the Gauss hypergeometric function:
$$
F(\alpha,\beta,\gamma;z)=\sum_{n\ge 0}z^n
\prod_{j=0}^{n-1}\frac{(\alpha+j)(\beta+j)}{(1+j)(\gamma+j)},
$$
and 
$$
\alpha=-1/2,\ \beta=2b,\ \gamma=2b+3/2.
$$
\end{proposition}

\begin{example} Consider the special case $b=1/2$ (the angle is
$\pi$, i.e. its boundary is the imaginary axis). In this case the
solution is 
$$
f(\zeta)=\zeta F(-1/2,1,5/2;-\zeta^2),
$$
i.e., 
$$
f(\zeta)=\frac{3}{8}
((\zeta-\zeta^{-1})+(\zeta+\zeta^{-1})^2{\rm arctan}\zeta).
$$
\end{example}

\begin{remark}
This analysis is similar to the analysis in \cite{GV}, subsection
2.2.3, where solutions are also expressed via the 
Gauss hypergeometric function.
\end{remark}

Lionel Levine has simulated on a computer the corresponding discrete (rotor-router) model 
for $b=1/2$ and $b=1/4$ (the dynamics in the half-plane and the quarter-plane, respectively),
and got a very close fit with formula (\ref{angl}). This gives rise to a conjecture that for these values of $b$ (and 
perhaps for all) the scaling limit is given by formula (\ref{angl}). 

\subsection{The killing and reflecting case}\label{kr}

Another problem one can consider is the 
same as the previous one, except that the angle is $\lbrace{z\in \Bbb C|{\rm arg}z\in
[0,\pi b]\rbrace}$ (i.e., the upper half of the angle considered before), and on the right side of the angle (the positive
half-axis) we have not the Dirichlet boundary condition for
$\Phi$, but the Neumann boundary condition, $\frac{\partial
\Phi}{\partial {\bold n}}=0$. This can be conjecturedly
interpreted as the scaling limit of the random walk problem on a
rectangular lattice in this angle with a killing condition on the
left side and reflecting condition on the right side. 

The problem of finding $\Omega$ in this case reduces to the previous
one, by doubling the angle using reflection with respect to the
positive half-axis. In particular, if $b=2$ (the angle is the
full plane, and we have the reflecting condition above the 
positive half-axis and the killing condition below), then the
conformal map of the disk with a cut at the positive half-axis to the
region $\Omega$ is 
$$ 
f(\zeta)=\zeta F(-1/2,4,11/2;-\zeta^{1/2}), 
$$ 
where $\zeta^{1/2}:=r^{1/2}e^{i\theta/2}$, $\theta\in [0,2\pi)$.

\section{Flows with the killing-passing-reflecting condition}

Consider now the discrete model as above, with the condition on the positive half-axis saying that 
particles coming from above are reflected, while particles coming from below pass through with probability $p$ and are killed with probability 
$1-p$. If $p=0$, this is the killing-reflecting case considered in the previous section, 
and if $p=1$, this is the passing-reflecting case, whose simulation is shown on the right side of Fig. 4 of \cite{LP}.  

We expect that the scaling limit of this process exists for any $p\in [0,1]$, 
and is described as follows. 

\begin{problem}\label{pr3} 
Let $\Omega$ be the asymptotic region. 
Then the boundary of $\Omega$ consists of a curve $\Gamma$ and the interval $[a,b]$ of the positive 
half-axis (the curve $\Gamma$ goes counterclockwise from $b$ to $a$ around the origin). Then the boundary conditions for the potential $\Phi$ in $\Omega$ 
are: 
$$
\Phi=0 \text{ on }\Gamma;
$$
$$
(\Phi_y)_+=0\text{ on }[a,b];
$$
$$
\Phi_-=0\text{ on }[0,a];
$$
$$
(\Phi_y)_+=p(\Phi_y)_-\text{ on }[0,a].
$$
Here $\Phi_y$ is the derivative of $\Phi$ with respect to $y$, and $+,-$ denote the one-sided limits 
from above and below ($\Phi$ and $\Phi_y$ are not assumed continuous on $(0,a)$).  

The condition on $\Omega$, as before, is the homothetic transformation under the flow defined by the normal derivative of $\Phi$. 
\end{problem}

\begin{remark} 1. To motivate these boundary conditions, it is convenient to use 
the following reformulation of the killing-passing-reflecting condition in the discrete model:
all the particles hitting the positive half-axis from below are killed, but for each $k$ particles killed 
on a small interval $I$ we produce about $pk$ particles which originate at $I$ and move in the upward direction. 
This explains the second and third boundary conditions 
(reflecting condition on $[a,b]$ from above, killing condition on $[0,a]$ from below),
and clarifies the meaning of the fourth boundary condition, which says that the flux through $[0,a]$ on the upper side is 
$p$ times the flux on the lower side. 

2. Note that for $p=0$, Problem \ref{pr3} reduces to the problem of Subsection \ref{kr} for $b=2$. 
\end{remark}

Unfortunately, we did not manage to compute the solution $\Omega$ of Problem \ref{pr3}
explicitly. However, here is an interesting moment property of this solution (which we expect to determine $\Omega$ uniquely up to scaling). 

Define a continuous branch of $\log(z)$ in $\Omega\setminus \Bbb R_+$ by defining it to be real on the reflecting boundary 
(i.e. the limit of $\log(z)$ from above is real for $z\in \Bbb R_+, z\ne 0$). This also defines branches of $z^s=e^{s\log(z)}$ for any $s$.  

\begin{proposition}\label{prop3}
(i) If $\Omega$ is a solution of Problem \ref{pr3}, then one has 
$$
\int_\Omega {\rm Re}(z^{n+\alpha})dxdy=0,\ \int_\Omega {\rm Re}(z^{n-\alpha})dxdy=0, \ n\in \Bbb N,
$$
where $\alpha={\rm arccos}(p)/2\pi$.
In particular, $0\le \alpha\le 1/4$.

(ii)
For $p=1$ ($\alpha=0$), one has  
$$
\int_\Omega {\rm Re}(z^n)dxdy=0,\  \int_\Omega {\rm Re}(z^n\log(z))dxdy=0,\ n\in \Bbb N.
$$
\end{proposition}

\begin{proof} (i) Let $u$ be a function on $\Omega\setminus \Bbb R_+$. 
Then the time derivative of the moment $\int_\Omega udxdy$ under the Hele-Shaw flow is given by 
$$
\frac{d}{dt}\int_\Omega udxdy=\int_{\Gamma}u\frac{\partial \Phi}{\partial \bold n}dl=
$$
$$
\int_{\partial \Omega}u\frac{\partial \Phi}{\partial \bold n}dl
-\int_0^bu_+(\Phi_y)_+dx+\int_0^au_-(\Phi_y)_-dx.
$$
By Green's formula, 
this equals
$$
\int_{\Omega} (u\Delta \Phi-\Phi\Delta u)dxdy+\int_{\partial \Omega}\Phi\frac{\partial u}{\partial \bold n}dl
-\int_0^bu_+(\Phi_y)_+dx+\int_0^au_-(\Phi_y)_-dx.
$$
If $u$ is harmonic and vanishes to sufficient order at $0$ then the first summand is zero. So, using the first three boundary conditions, we have 
\begin{equation}\label{timeder}
\frac{d}{dt}\int_\Omega udxdy=\int_0^b\Phi_+(u_y)_+dx
-\int_0^au_+(\Phi_y)_+dx+\int_0^au_-(\Phi_y)_-dx.
\end{equation}
Now suppose that $u={\rm Re}(z^{n\pm \alpha})=r^{n+\alpha}{\rm cos}((n\pm\alpha)\theta)$, $z=re^{i\theta}$, $\theta\in [0,2\pi)$. 
Then $(u_y)_+=0$, so the first summand in (\ref{timeder}) vanishes. Thus, we have 
\begin{equation}\label{timeder1}
\frac{d}{dt}\int_\Omega udxdy=\int_0^au_-(\Phi_y)_-dx
-\int_0^au_+(\Phi_y)_+dx+.
\end{equation}
But we have $u_+=x^\alpha$, $u_-=x^\alpha {\rm cos}(2\pi \alpha)=px^\alpha$. 
So by the fourth boundary condition, we get from (\ref{timeder}) that 
$\frac{d}{dt}\int_\Omega udxdy=0$. But since the time evolution of $\Omega$ is homothetic, and $\int_\Omega udxdy$ is homogeneous
of nonzero degree under dilations, we get $\int_\Omega udxdy=0$, as desired.   

(ii) follows from (i) by considering the difference of the two equations in (i) 
and computing the first term of the Taylor expansion in $\alpha$ at $\alpha=0$. 
\end{proof}

\begin{remark}
Let $\Omega$ be a region solving Problem \ref{pr3} for $p=1$ (which is expected to match the right shape 
at Fig. 4 of \cite{LP}), and $\overline{\Omega}$ be the complex conjugate region. 
Proposition \ref{prop3}(ii) implies that 
$$
\int_\Omega z^ndxdy+\int_{\overline\Omega}z^ndxdy=0,\ n\in \Bbb N.
$$
This implies that the ``balayage sum'' $E$ of $\Omega$ and $\overline{\Omega}$ (i.e., 
the region whose moments are equal to the sum of the moments of $\Omega$ 
and of $\overline{\Omega}$, see \cite{GV}) satisfies the condition 
$\int_E z^ndxdy=0$, $n\in \Bbb N$, i.e. is a disk (see \cite{VE} and references therein). 
\end{remark}

Note that for $p=0$ ($\alpha=1/4$), i.e. the killing-reflecting case discussed in Subsection \ref{kr}, 
the two arithmetic progressions in Proposition \ref{prop3}(i) combine into a single one, and yield
$$
\int_\Omega {\rm Re}(z^{n/2+1/4})dxdy=0, \ n\in \Bbb N.
$$
This can be generalized to the case of any angle (i.e., any $b$). Namely, 
we have the following proposition.

\begin{proposition}
Let $\Omega$ be a region solving Problem \ref{pr2}. Then 
$$
\int_\Omega z^{\frac{2n+1}{2b}}dxdy=0, n\in \Bbb N.  
$$
Therefore, the solution of the killing-reflecting problem with parameter $b$
(which is the upper half of the region $\Omega$) satisfies the equality
$$
\int_\Omega {\rm Re}(z^{\frac{2n+1}{2b}})dxdy=0, n\in \Bbb N.  
$$
\end{proposition}

The proof of this proposition is similar to the proof of Proposition \ref{prop3}, using Green's formula. 

A similar property was established by Lionel Levine for the discrete model (for $b=1$).

\end{document}